\documentclass[11pt]{article}

\usepackage{ccicons}
\usepackage{amssymb}
\usepackage{amsmath}
\usepackage[dvips]{epsfig}
\usepackage{caption2}
\usepackage{graphicx}
\usepackage[colorlinks=true,linkcolor=blue,citecolor=red]{hyperref}

\newtheorem{Lemma}{Lemma}[section]
\newtheorem{Theorem}{Theorem}
\newtheorem{Proposition}[Lemma]{Proposition}

 {\begin{trivlist} \item[]{\bf Proof. }}%
 {\hspace*{\fill}$\rule{.4\baselineskip}{.4\baselineskip}$\end{trivlist}}

 {\begin{trivlist}\item[]\textbf{Acknowledgments }}{\end{trivlist}}

\makeatletter\@addtoreset{figure}{section}\makeatother

\makeatletter \@addtoreset{equation}{section} \makeatother

\newcommand{\R}{\mathbb{R}}
\newcommand{\C}{\mathbb{C}}

\newcommand{\Z}{\mathbb{Z}}

\def\Re{\mathop{\mathrm{Re}}}
\def\Im{\mathop{\mathrm{Im}}}

\newcommand{\rme}{\mathrm{e}}
\newcommand{\rmi}{\mathrm{i}}

\renewcommand{\leq}{\leqslant}
\renewcommand{\geq}{\geqslant}

\newcommand{\ess}{\mathrm{ess}}
\newcommand{\per}{\mathrm{per}}
\newcommand{\loc}{\mathrm{loc}}

\newfam\bifam
\font\tenbi=cmmib10 scaled \magstep1 \font\sevenbi=cmmib10 at 11pt
\font\fivebi=cmmib10 at 6pt \textfont\bifam = \tenbi
\scriptfont\bifam = \sevenbi \scriptscriptfont\bifam= \fivebi

\begin{document}

\thispagestyle{empty}

\title{\bf Transverse instability of periodic and generalized solitary waves for a fifth-order KP model} 

\author{\Large Mariana Haragus\footnote{Institut FEMTO-ST and LMB, Univ. Bourgogne Franche-Comt\'e, 25030 Besan\c{c}on cedex, France} 
\ \&\ 
Erik Wahl\'en\footnote{Centre for Mathematical Sciences, Lund University, P.O. Box 118, 22100 Lund, Sweden} }

\date{} 

\maketitle 

\begin{abstract}
We consider a fifth-order Kadomtsev-Petviashvili equation which arises as a two-dimensional model in the classical water-wave problem. This equation possesses a family of generalized line solitary waves which decay exponentially to periodic waves at infinity. We prove that these solitary waves are transversely spectrally unstable and that this instability is induced by the transverse instability of the periodic tails. We rely upon a detailed spectral analysis of some suitably chosen linear operators.
\end{abstract}

{\small

{\bf Keywords:} {transverse stability, periodic waves, generalized solitary waves, dispersive equations} 

Accepted for publication in Journal of Differential Equations. DOI:\href{http://dx.doi.org/10.1016/j.jde.2016.11.025}{10.1016/j.jde.2016.11.025}

This work is licensed under a CC BY-NC-ND 4.0 license.
}

\section{Introduction}

We consider a fifth-order Kadomtsev-Petviashvili (KP) equation
\begin{equation}\label{e:kp6}
\partial_t\partial_x u = \partial_x^2 \left(
\partial_x^4 u + \partial_x^2 u + \frac12 u^2
\right) + \partial_y^2 u,
\end{equation}
in which the unknown $u$ depends upon two space variables $(x,y)\in\R^2$ and time $t\in\R$. This equation arises as a two-dimensional model for capillary-gravity water waves in the regime of critical surface tension, when the Bond number is close to $1/3$ \cite{La}. While the exact values of the coefficients in \eqref{e:kp6} are unimportant, their signs have been chosen corresponding to the case of Bond number less than $1/3$. We can also regard \eqref{e:kp6} as a two-dimensional version of the Kawahara equation
\begin{equation}\label{e:kaw}
\partial_tu = \partial_x \left(
\partial_x^4 u + \partial_x^2 u + \frac12 u^2
\right),
\end{equation}
just as the KP equation is a two-dimensional version of the well-known Korteweg-de Vries (KdV) equation.

The Kawahara equation \eqref{e:kaw} possesses a family of traveling generalized solitary waves \cite{HS}. In contrast to the solitary waves of the KdV equation which tend exponentially to zero as $|x|\to\infty$, the generalized solitary waves of \eqref{e:kaw} decay exponentially to periodic waves as  $|x|\to\infty$. The amplitude of these periodic waves may be exponentially small but not zero \cite{Lo}. Our purpose is to study the transverse stability of these generalized solitary waves, i.e., their stability as solutions of the equation \eqref{e:kp6}, for perturbations which depend upon both spatial variables $x$ and $y$. While the stability of solitary waves for the KdV equation has been intensively studied, very little is known about the stability of the generalized solitary waves for \eqref{e:kaw}. A key difference is that asymptotically these generalized solitary waves tend to periodic waves, not to a constant, and the stability of these periodic waves is not fully understood. We mention the result in \cite{HLS} showing that these periodic waves are spectrally stable as solutions of \eqref{e:kaw} provided their amplitude is sufficiently small. 

Our main result shows that the generalized solitary waves of \eqref{e:kaw} are transversely spectrally unstable. The starting point of our analysis is a formulation of the transverse instability problem in terms of the spectrum of a suitably chosen operator. In particular, this allows to distinguish between linear instability, due to point spectrum, and essential instability, due to essential spectrum. Next, the key step is the spectral analysis of the operator found for the asymptotic periodic waves. We prove that these periodic waves are transversely linearly unstable with respect to perturbations which are co-periodic (i.e., they have the same period as the periodic wave) in the longitudinal direction, and transversely essentially unstable with respect to perturbations which are localized in the longitudinal direction. Finally, a rather general perturbation argument allows to conclude to the essential instability of the generalized solitary waves.

It is interesting to compare these properties with known stability results for the KP equation. Recall that the KP equation comes in two flavors:  KP-I, which is valid for strong surface tension (Bond number greater than $1/3$), and KP-II, which is valid for weak or zero surface tension (Bond number less than $1/3$) \cite{La}. Both equations reduce to the KdV equation in  one dimension, but display completely opposite transverse dynamics. While KP-I predicts transverse instability of both periodic and solitary waves, KP-II predicts stability \cite{APS, hss, Ha2, jz, KP, KSF, M, MT, RT1}. The latter may at first sight seem to contradict our result, but note that the KP-II equation does not capture the small periodic tails of the generalized solitary waves (nor is it uniformly valid in the limit of critical surface tension).

Generalized line solitary waves with small amplitude are known to exist for the full capillary-gravity water-wave problem in the regime of weak surface tension \cite{IK, Lo}. When the surface tension is close to critical, these are to leading order described by the Kawahara equation \eqref{e:kaw}. The results in this paper predict that these generalized solitary waves are transversely unstable and that this instability is due to that of the asymptotic periodic waves. This prediction will make the object of future work. We point out that the instability predictions based on the KP-I equation for the regime of strong surface tension have been confirmed for the full water-wave problem for both periodic and solitary waves \cite{GHS, Ha3, PS}. Also notice that in the regime of weak surface tension, transverse instability of solitary waves have been recently proved for a class of true solitary waves which appear in the modulational regime, which is different from the long-wave setting studied here \cite{GSW}.

Finally, the question of nonlinear transverse instability for these generalized solitary waves remains open. This question has recently been solved in the positive for line solitary waves of the KP-I equation \cite{RT1}, as well as the full water wave problem with strong surface tension \cite{RT2}. Although it would seem natural to expect a similar result for the fifth-order KP equation, the proof is far from straightforward for generalized solitary waves, since the instability is due to essential spectrum (even considering solutions which are periodic in the transverse direction). 

The paper is organized as follows. In Section 2 we recall the existence results for both periodic and generalized solitary waves of the Kawahara equation \eqref{e:kaw}. The main results are presented in Section 3, and the proofs are given in Section 4.

\section{Existence of line traveling waves}
\label{s:exist}

In this section we recall the existence results for both periodic and generalized solitary traveling waves of the equation \eqref{e:kaw}. 

We consider traveling waves moving with constant speed $c$. In a comoving frame, after replacing $x+ct$ by $x$, these traveling waves are stationary solutions of the equation
\begin{equation}\label{e:star}
\partial_tu = \partial_x \left(
\partial_x^4 u + \partial_x^2 u - cu + \frac12 u^2
\right),
\end{equation}
and therefore satisfy the ODE
\[
\partial_x^4 u + \partial_x^2 u - cu + \frac12 u^2 = C,
\]
obtained after integrating the right hand side of \eqref{e:star} once with respect to $x$. As a remnant of the Galilean invariance of \eqref{e:kaw}, we may set $C=0$ and restrict to the solutions of the ODE
\begin{equation}\label{e:2}
\partial_x^4 u + \partial_x^2 u - cu + \frac12 u^2 = 0.
\end{equation}

\subsection{Small periodic waves}

The existence of small periodic solutions of \eqref{e:2} for small speeds $c$ has been proved in \cite{HLS} (see also \cite[Chapters 4 and 7]{Lo}). We recall in the next proposition the result from \cite{HLS} which also gives a number of properties of these periodic solutions which are essential in our analysis.

\begin{Proposition}{\bf (\cite[Theorem 1]{HLS})} \label{p:1}
There exist positive constants $c_0$ and $a_0$ such that, for any $c\in(-c_0,c_0)$, the equation \eqref{e:2} possesses a one-parameter family of even, periodic solutions $(\varphi_{a,c})_{a\in(-a_0,a_0)}$ of the form
\[
\varphi_{a,c}(x) = p_{a,c}(k_{a,c}x),\quad \forall\ x\in\R,
\]
with the following properties.
\begin{enumerate}
\item The real-valued map $(a,c)\mapsto k_{a,c}$ is analytic on $(-a_0,a_0)\times(-c_0,c_0)$ and
\[
k_{a,c} = k_0(c)+c\widetilde k(a,c),\quad 
k_0(c) = \left(\frac{1+\sqrt{1+4c}}2\right)^{1/2},\quad
\widetilde k(a,c) = \sum_{n\geq1}\widetilde k_{2n}(c)a^{2n},
\]
for any $(a,c)\in(-a_0,a_0)\times(-c_0,c_0)$, where $|\widetilde k_{2n}(c)|\leq K_0/\rho_0^{2n}$, for any $n\geq1$ and some positive constants $K_0$ and $\rho_0$.
\item The map $(a,c)\mapsto p_{a,c}$ is analytic on $(-a_0,a_0)\times(-c_0,c_0)$ with values $H^6_{\per}(0,2\pi)$ and 
\[
p_{a,c}(z) = ac \cos(z) + c \sum_{\substack{n,m\geq0,n+m\geq2\\n-m\not=\pm1}}
\widetilde p_{n,m}(c)\rme^{\rmi(n-m)z}a^{n+m},
\]
in which $\widetilde p_{n,m}(c)$ are real numbers such that $\widetilde p_{n,m}(c) = \widetilde p_{m,n}(c)$ and $|\widetilde p_{n,m}(c)|\leq C_0/\rho_0^{n+m}$, for any $c\in(-c_0,c_0)$ and some positive constant $C_0$. (Here $H^{6}_{\mathrm{per}}(0,2\pi)$ is the space of $2\pi$-periodic functions defined in Section \ref{ss:periodic waves} below.)
\item The Fourier coefficients $\widehat p_q(a,c)$ of the $2\pi$-periodic function $p_{a,c}$,
\[
p_{a,c}(z) = \sum_{q\in\Z} \widehat p_q(a,c)\rme^{\rmi qz},\quad
\forall\  z\in\R,
\]
are real and satisfy $\widehat p_0(a,c)=O(ca^2)$ and $\widehat p_q(a,c) = O(c|a|^{|q|})$, for all $q\not=0$, as $a\to0$. Moreover, the map $a\mapsto \widehat p_q(a,c)$ is even (resp.~odd) for even (resp.~odd) values of $q$, and in particular $p_{-a,c}(z) = p_{a,c}(z+\pi)$.
\end{enumerate}
\end{Proposition}

We collect below some properties of  $k_{a,c}$ and $p_{a,c}$ which are needed in our proofs. First, a direct calculation allows to compute the expansions of $k_{a,c}$ and $p_{a,c}$, as $a\to0$. Without writing explicitly the dependence upon $c$, for notational simplicity, we find
\begin{equation}\label{e:3}
k_{a,c}=k_0+\sum_{n\geq1}k_{2n}a^{2n},\quad
k_0^2 = \frac{1+\sqrt{1+4c}}2,\quad
k_2(4k_0^3-2k_0) = -\frac c4 + \frac{c^2}{8X_2},
\end{equation}
in which we used the notation
\begin{equation}\label{e:xn}
X_n = k_0^4n^4 -k_0^2n^2 -c,\quad \forall\ n\geq 2.
\end{equation}
For $p_{a,c}$, we write
\begin{equation}\label{e:4}
p_{a,c}(z) = ac\left( \cos(z) + \sum_{n\geq1} p_n(z)a^n\right),
\end{equation}
in which we find
\[
p_1(z) = \frac14 - \frac c{4X_2}\cos(2z),\quad
p_2(z) = \frac{c^2}{8X_2X_3}\cos(3z).
\]
We point out that the Fourier coefficients $\pm1$ of the functions $p_n(z)$ are zero by construction. 

Next, notice that $k_{a,0}=1$, $p_{a,0}=0$, and we claim that
\begin{equation}\label{e:e5}
\partial_ck_{a,c}^2|_{c=0} = 1-q(a), \quad
\partial_cp_{a,c}|_{c=0} = a \cos(z) + q(a),\quad q(a)=1-\sqrt{1-\frac12 a^2}.
\end{equation}
Indeed, recall that
\begin{equation}\label{e:e4}
k_{a,c}^4\partial_z^4p_{a,c} + k_{a,c}^2\partial_z^2p_{a,c} - cp_{a,c} + \frac12 p_{a,c}^2 = 0.
\end{equation}
Differentiating this equality with respect to $c$ and taking $c=0$ we find
\[
\left( \partial_z^4+\partial_z^2\right)\left(\partial_cp_{a,c}|_{c=0}\right) =0,
\]
so that $\partial_cp_{a,c}|_{c=0}$ belongs to the kernel of $ \partial_z^4+\partial_z^2$. Since  $\partial_cp_{a,c}|_{c=0}$ is an even function, this implies that  $\partial_cp_{a,c}|_{c=0}$ is a linear combination of $\cos(z)$ and $1$, and taking into account the expansion in Proposition~\ref{p:1} (ii), we obtain the second equality in \eqref{e:e5}. Next, we differentiate \eqref{e:e4} twice with respect to $c$ and take $c=0$. This gives 
\[
\left( \partial_z^4+\partial_z^2\right)\left(\partial_c^2p_{a,c}|_{c=0}\right)
+2a\left( \partial_ck_{a,c}^2|_{c=0} 
 -1
+q(a)\right)\cos(z)+q(a)^2-2q(a)+a^2\cos^2(z)=0,
\]
and the solvability conditions for this equation imply the first and the third equalities in \eqref{e:e5}.

\subsection{Generalized solitary waves}

The existence of generalized solitary waves is a consequence of the result in \cite[Chapter 7, Theorem 7.1.18]{Lo} for general four-dimensional reversible ODEs in presence of a $0^2(\rmi\omega)$ resonance. For completeness, we give the proof of the following proposition in Appendix~\ref{a:exist_gw}.

\begin{Proposition} \label{p:2}
There exist positive constants $a_1$ and $M_1$ such that for any $0<\ell<\pi$ and $0<\lambda<1$, there exist $c_2(\ell)>0$ and $a_2(\ell)>0$ such that for all $c\in(0,c_2(\ell)]$ and $|a|\in[a_2(\ell)c\rme^{-\ell/\sqrt c},a_1]$, the equation \eqref{e:2} possesses an even solution
\begin{equation}\label{e:6bis}
u_{a,c}(x) = h_{a,c}(x)+\varphi_{a,c}(x+\tau_{a,c}\tanh(\sqrt c x/2)),
\end{equation}
with the following properties:
\begin{enumerate}
\item $|\partial_x^j h_{a,c}(x)|\leq M_1 c \rme^{-\lambda\sqrt c |x|}$ for $j=0,1,2,3$ and all $x\in\R$;
\item $\varphi_{a,c}$ is the periodic solution in Proposition~\ref{p:1};
\item the asymptotic phase shift is such that $\tau_{a,c} = O(1)$, as $(a,c)\to(0,0)$.
\end{enumerate}
\end{Proposition}

\section{Transverse instability: main results}
\label{s:main}

In this section, we state the main instability results for both periodic and generalized solitary waves. We give the proofs of these results in Section~\ref{s:proofs}.

\subsection{Formulation of the transverse instability problem}

Assume that $u_*$ is a one-dimensional solution of \eqref{e:star}, for instance a periodic wave (as in Proposition~\ref{p:1}) or a generalized solitary waves (as in Proposition~\ref{p:2}). Consider the linearized equation
\begin{equation}\label{e:7}
\partial_t\partial_x u =
\partial_x^2\left(\partial_x^4u + \partial_x^2u -cu + u_*u\right)
+\partial_y^2u,
\end{equation}
and set
\[
\mathcal A_* = \partial_x^2\left(\partial_x^4 + \partial_x^2 -c + u_*\right).
\]
Roughly speaking, the wave $u_*$ is called transversely unstable if the equation \eqref{e:7} possesses solutions of the form
\[
u(t,x,y) = \rme^{\lambda t}v(x,y),
\]
for some $\Re\lambda>0$ and $v$ a time-independent function which belongs to the set of the allowed perturbations. Since $u_*$ does not depend upon the transverse spatial variable $y$, the operator $\mathcal A_*+\partial_y^2$ in the right hand side of \eqref{e:7} has $y$-independent coefficients, so that using Fourier transform in $y$ we can reformulate the instability statement and say the $u_*$ is {transversely unstable} if the linearized equation
\[
\partial_t\partial_xu = \mathcal A_*u-\omega^2u,
\]
has solutions of the form
\[
u(t,x) = \rme^{\lambda t}v(x),
\]
for some $\Re\lambda>0$, $\omega\in\R$, and $v$ in some space $H$ of functions depending upon the longitudinal spatial variable $x$, only. In this setting, perturbations are bounded in the transverse variable $y$ and determined by the choice of $H$ in the longitudinal variable $x$ (e.g., localized if $H=L^2(\R)$ or periodic if $H=L^2(0,L)$). 

We can now reformulate the transverse instability problem and say that \emph{$u_*$ is transversely spectrally unstable} if the linear operator
$\lambda \partial_x -\mathcal A_* + \omega^2$,
is not invertible in $H$ for some $\lambda\in\C$ with $\Re\lambda>0$ and $\omega\in\R$. The particular form of this operator allows to further say that $u_*$ is transversely spectrally unstable if the spectrum of the linear operator $\lambda\partial_x-\mathcal A_*$ contains a negative value $-\omega^2$ for some $\Re\lambda>0$. Notice that if $-\omega^2$ is an isolated eigenvalue of $\lambda\partial_x-\mathcal A_*$ then this definition implies transverse linear instability. If $-\omega^2$ belongs to the essential spectrum of $\lambda\partial_x-\mathcal A_*$, 
\begin{equation}\label{e:8}
\sigma_{\ess}(\lambda\partial_x-\mathcal A_*) = 
\{ \nu\in\C\;;\; \lambda\partial_x-\mathcal A_*-\nu 
\mbox{ is not Fredholm with index } 0\},
\end{equation}
we may say that $u_*$ is transversely essentially unstable.
Our main results show that the periodic waves $\varphi_{a,c}$ are transversely unstable with respect to co-periodic longitudinal perturbations (Theorem~\ref{t:1}), and that both the periodic waves $\varphi_{a,c}$ and the generalized solitary waves $u_{a,c}$ are transversely essentially unstable with respect to localized longitudinal perturbations (Theorem~\ref{t:2} and Theorem~\ref{t:3}, respectively).

\subsection{Periodic waves}
\label{ss:periodic waves}

Consider the small periodic waves $\varphi_{a,c}$ constructed in Proposition~\ref{p:1}, and the linear operator
\[
\lambda\partial_x-\mathcal A_{a,c},\quad
\mathcal A_{a,c} = \partial_x^2\left(\partial_x^4+\partial_x^2-c+\varphi_{a,c}\right).
\]
Since the period $2\pi/k_{a,c}$ of $\varphi_{a,c}$ depends upon $a$ and $c$, it is convenient to rescale $x$ and $\lambda$ by taking
\[
z=k_{a,c}x,\quad \lambda = k_{a,c}\Lambda,
\]
and work with the rescaled operator
\[
\Lambda\partial_z-\mathcal B_{a,c},\quad
\mathcal B_{a,c} = \partial_z^2(k_{a,c}^4\partial_z^4+k_{a,c}^2\partial_z^2-c+p_{a,c}),
\]
which has $2\pi$-periodic coefficients. 

First, for co-periodic perturbations we take $H=L^2(0,2\pi)$ and as domain of definition for $\mathcal B_{a,c}$ the subspace $H^6_{\per}(0,2\pi)$ consisting of $2\pi$-periodic functions,
\[
H^j_{\per}(0,2\pi) = \{ f\in H^j_\loc(\R) \;;\; f(z+2\pi) = f(z),\ \forall\ z\in\R\},
\]
for $j\geq1$. Then $\mathcal B_{a,c}$ is closed in $H$, and our main result is the following theorem which is proved in Section~\ref{ss:pt1}.  

\begin{Theorem}\label{t:1}
There exist positive constants $c_3$ and $a_3$ such that for any $c\in(-c_3,c_3)$ and $a\in(-a_3,a_3)$, there exists $\Lambda_{a,c}>0$ such that for any $\Lambda\in(0,\Lambda_{a,c})$ the linear operator $\Lambda\partial_z-\mathcal B_{a,c}$ acting in $L^2(0,2\pi)$ with domain  $H^6_{\per}(0,2\pi)$ has a simple negative eigenvalue. Consequently, the periodic wave $\varphi_{a,c}$ is transversely linearly unstable with respect to co-periodic longitudinal perturbations.
\end{Theorem}

Next, for localized perturbations we take $H=L^2(\R)$ and as domain of definition of $\mathcal B_{a,c}$, and also $\Lambda\partial_z-\mathcal B_{a,c}$, the subspace $H^6(\R)$. The spectral analysis in this space is based on a Bloch-wave (or Floquet in this case) decomposition, which shows that the spectrum of $\Lambda\partial_z-\mathcal B_{a,c}$ in $L^2(\R)$ is the union of the spectra of the operators
\[
\Lambda(\partial_z+\rmi\gamma)-\mathcal B_{a,c,\gamma},\quad
\mathcal B_{a,c,\gamma} = (\partial_z+\rmi\gamma)^2 \left(
k_{a,c}^4(\partial_z+\rmi\gamma)^4 + k_{a,c}^2(\partial_z+\rmi\gamma)^2
-c+p_{a,c}\right),
\]
acting in $L^2(0,2\pi)$ with domain $H^6_\per(0,2\pi)$, for $\gamma\in(-1/2,1/2]$ (e.g., see \cite{Ha1, Jo}). Then the transverse spectral instability of $\varphi_{a,c}$ with respect to localized perturbations is an immediate consequence of Theorem~\ref{t:1}. Moreover, from \cite{Ha1, Jo} we deduce that the spectrum is purely essential spectrum, in the sense of definition \eqref{e:8}, so that the instability is essential. Summarizing, we have the result below.

\begin{Theorem}\label{t:2}
For any $a$, $c$, and $\Lambda$ as in Theorem~\ref{t:1}, the linear operator $\Lambda\partial_z-\mathcal B_{a,c}$ acting in $L^2(\R)$ with domain $H^6(\R)$ has negative essential spectrum. Consequently, the periodic wave $\varphi_{a,c}$ is transversely essentially unstable with respect to localized longitudinal perturbations.
\end{Theorem}

\subsection{Generalized solitary waves}

Consider the generalized solitary waves found in Proposition~\ref{p:2} and the linear operator
\[
\lambda\partial_x-\mathcal C_{a,c},\quad
\mathcal C_{a,c} = \partial_x^2\left(\partial_x^4+\partial_x^2-c+u_{a,c}\right),
\]
acting in $L^2(\R)$ with domain $H^6(\R)$ (localized perturbations). The key observation in the spectral analysis of $\lambda\partial_x-\mathcal C_{a,c}$ is that it is a relatively compact perturbation of the asymptotic operator
\[
\lambda\partial_x-\mathcal C_{a,c}^\infty = \left\{ \begin{array}{ll}
\lambda\partial_x-\mathcal A_{a,c}^+, & \mbox{ for } x>0\\
\lambda\partial_x-\mathcal A_{a,c}^-, & \mbox{ for } x<0
\end{array}\right. ,\qquad
\mathcal A_{a,c}^\pm = \partial_x^2\left(\partial_x^4+\partial_x^2 -c +
\varphi_{a,c}(\cdot\pm\tau_{a,c})\right).
\]
Since the essential spectrum is stable under relatively compact perturbations, the generalized solitary wave $\varphi_{a,c}$ is transversely unstable provided the asymptotic operator has negative essential spectrum. The latter property is a consequence of Theorem~\ref{t:1}, just as Theorem~\ref{t:2}. This leads to the following result, which is proved in Section~\ref{ss:pt3} 

\begin{Theorem}\label{t:3}
For any $a$, $c$, and $\Lambda$ as in Proposition~\ref{p:2} and Theorem~\ref{t:1}, the linear operator $k_{a,c}\Lambda\partial_x-\mathcal C_{a,c}$ acting in $L^2(\R)$ with domain $H^6(\R)$ has negative essential spectrum. Consequently, the generalized solitary wave $u_{a,c}$ is transversely essentially unstable.
\end{Theorem}

\section{Proofs}
\label{s:proofs}

\subsection{Proof of Theorem~\ref{t:1}}
\label{ss:pt1}

We claim that it is enough to prove that the operator  $\mathcal B_{a,c}$ has a simple positive eigenvalue, or equivalently, that the operator  $-\mathcal B_{a,c}$ has a simple negative eigenvalue, for sufficiently small $a$ and $c$. Indeed, notice that the operator $\Lambda\partial_z-\mathcal B_{a,c}$ is a small relatively bounded perturbation of $-\mathcal B_{a,c}$, for any sufficiently small $\Lambda\in\C$. If $-\mathcal B_{a,c}$ has a simple negative eigenvalue, then a standard perturbation argument implies that $\Lambda\partial_z-\mathcal B_{a,c}$ has a simple eigenvalue in some open disk centered on the real axis and contained in the open left half complex plane. For real values $\Lambda$, the operator $\Lambda\partial_z-\mathcal B_{a,c}$ is real, so its spectrum is symmetric with respect to the real axis. Consequently, the simple eigenvalue above is necessarily real and negative, which proves the claim. 

For small $a$ and $c$, the operator $\mathcal B_{a,c}$ is a small relatively bounded perturbation of the operator
\[
\mathcal B_{0,0} = \partial_z^2\left(\partial_z^4 + \partial_z^2\right),
\]
which has constant coefficients. We can compute the spectrum of $\mathcal B_{0,0}$ using Fourier series, and find 
\[
\sigma(\mathcal B_{0,0}) = \{ -n^2(n^4-n^2),\ n\in\Z\},
\]
where $0$ is a semi-simple triple eigenvalue and all other eigenvalues are negative.  Then a standard perturbation argument shows that there exists a neighborhood $V$ of $0$ in the complex plane and a positive constant $m$ such that 
$V\subset\{\nu\in\C\; ;\; |\Re\nu|<m/2\}$ and
for sufficiently small $a$ and $c$, the spectrum of $\mathcal B_{a,c}$ decomposes as
\[
\sigma(\mathcal B_{a,c}) = \sigma_1(\mathcal B_{a,c}) \cup
\sigma_2(\mathcal B_{a,c}),\quad \sigma_1(\mathcal B_{a,c})\subset V,\quad
\sigma_2(\mathcal B_{a,c})\subset\{\nu\in\C\; ;\; \Re\nu<-m\},
\]
and $\sigma_1(\mathcal B_{a,c})$ contains precisely three eigenvalues, not necessarily distinct, counted with multiplicities. We show that one of these eigenvalues is positive when $a\not=0$ and that the other two eigenvalues are equal to $0$.

For $a=0$, the operator $\mathcal B_{0,c}$ has constant coefficients and using Fourier series, again, we can compute its spectrum,
\[
\sigma(\mathcal B_{0,c}) = \{ -n^2(k_0^2n^4-k_0^2n^2-c),\ n\in\Z\},
\]
where $k_0$ is the constant in the expansion \eqref{e:3} of $k_{a,c}$. In particular, $0$ is a triple eigenvalue of $\mathcal B_{0,c}$ with associated eigenfunctions $1$, $\cos(z)$, and $\sin(z)$. 

For $a\not=0$, we write
\[
\mathcal B_{a,c} = \partial_z^2\mathcal L_{a,c},\quad 
\mathcal L_{a,c}= k_{a,c}^4\partial_z^4 + k_{a,c}^2\partial_z^2 -c + p_{a,c}.
\]
The spectrum of $\partial_z\mathcal L_{a,c}$ has been studied in \cite{HLS}. According to \cite[Remark 3.5 (ii)]{HLS}, the kernel of $\partial_z\mathcal L_{a,c}$ is two-dimensional, spanned by the odd function $\partial_zp_{a,c}$ and an even function $\xi_{a,c}^e$ which is equal to $1$ when $a=0$. Since the kernel of $\mathcal B_{a,c}$ contains the kernel of $\partial_z\mathcal L_{a,c}$, $0$ is at least a double eigenvalue of $\mathcal B_{a,c}$ with two associated eigenfunctions $\xi_{a,c}^o$ and $\xi_{a,c}^e$ which are smooth continuations, for small $a$, of the vectors $\sin(z)$ and $1$, and are odd and even functions, respectively. 

In order to compute the third eigenvalue in $\sigma_1(\mathcal B_{a,c})$, we consider a basis for the associated three-dimensional spectral subspace which is a smooth continuation of the basis $\{1,\cos(z),\sin(z)\}$ found for $a=0$. Since $\mathcal B_{a,c}$ leaves invariant the subspaces consisting of even and odd functions, two vectors in this basis are even functions and the third one is an odd function. Clearly, the two eigenfunctions $\xi_{a,c}^o$ and $\xi_{a,c}^e$  above belong to this basis, and a third vector is an even function that we denote by $\psi_{a,c}$. Since $\xi_{a,c}^e=1+O(|a|)$ and $\psi_{a,c}=\cos(z)+O(|a|)$, upon replacing  $\psi_{a,c}$ by a linear combination of  $\psi_{a,c}$ and  $\xi_{a,c}^e$,  we can always choose $\psi_{a,c}$ to be orthogonal to $1$. Then, writing
\[
\mathcal B_{a,c}\psi_{a,c} = \nu_{a,c}\psi_{a,c} + \mu_{a,c}\xi_{a,c}^e,
\]
for small $a$, and taking the scalar product with $1$ we conclude that $\mu_{a,c}=0$. This implies that we can determine $\psi_{a,c}$ and the third eigenvalue $\nu_{a,c}$ by solving the eigenvalue problem
\begin{equation}\label{e:9}
\mathcal B_{a,c}\psi_{a,c} = \nu_{a,c}\psi_{a,c},
\end{equation}
in which $\nu_{a,c}$ and $\psi_{a,c}$ depend smoothly upon $a$ and $c$, and 
\[
\nu_{a,c}=O(|a|),\quad \psi_{a,c} = \cos(z)+O(|a|). 
\]
To complete the proof, we show that  $\nu_{a,c}$ is positive. 

First, recall that $p_{-a,c}(x) = p_{a,c}(x+\pi)$ which implies that $\nu_{-a,c}=\nu_{a,c}$, and in particular $\nu_{a,c}=O(a^2)$, as $a\to0$. Next, we claim that $\nu_{a,c}=O(c^2)$, as $c\to0$, so that $\nu_{a,c}=O(a^2c^2)$. Since $\mathcal B_{a,0} = \mathcal B_{0,0}$, we have that 
\begin{equation}\label{e:e2}
\nu_{a,0}=0,\quad \psi_{a,0} = \cos(z).
\end{equation}
Differentiating \eqref{e:9} with respect to $c$ and taking $c=0$ we find
\[
\mathcal B_{0,0}\left(\partial_c\psi_{a,c}|_{c=0}\right) + 
\partial_z^2\left( 2\partial_ck_{a,c}^2|_{c=0} \partial_z^4 + 
\partial_ck_{a,c}^2|_{c=0}\partial_z^2 - 1 + \partial_cp_{a,c}|_{c=0}\right)
\cos(z) = \partial_c\nu_{a,c}|_{c=0}\cos(z).
\]
The solvability condition for this equation gives
\[
\partial_c\nu_{a,c}|_{c=0} = - \partial_ck_{a,c}^2|_{c=0} 
+1 - \left[\partial_cp_{a,c}|_{c=0}\cos(z)\right]_1,
\]
in which the bracket $[u]_1$ represents the coefficient of $\cos(z)$ in the Fourier expansion of $u$. Taking into account the equalities \eqref{e:e5}, we conclude that $\partial_c\nu_{a,c}|_{c=0}=0$ which proves the claim. 

Next, consider the expansions
\[
\nu_{a,c}= \nu_2 a^2 + O(a^4),\quad \psi_{a,c}(z) = \cos(z) + \psi_1(z)a+\psi_2(z)a^2 + O(a^3),
\]
for small $a$. Inserting these expansions into \eqref{e:9}, at order $O(1)$ we find the eigenvalue problem at $a=0$, which holds, and at order $O(a)$ the equality
\[
\partial_z^2\mathcal L_0\psi_1 + \partial_z^2\mathcal L_1\cos(z) = 0,
\]
in which
\[
\mathcal L_0 = k_0^4\partial_z^4+k_0^2\partial_z^2-c,\quad
\mathcal L_1 = c\cos(z).
\]
A direct calculation gives
\[
\psi_1(z) =- \frac c{2X_2}\cos(2z),
\]
with $X_2$ given by \eqref{e:xn}. At order $O(a^2)$ we obtain
\[
\partial_z^2\mathcal L_0\psi_2 + \partial_z^2\mathcal L_1\psi_1 + \partial_z^2\mathcal L_2\cos(z) = \nu_2\cos(z),
\]
in which
\[
\mathcal L_2 = c\left(\frac14 - \frac c{4X_2}\cos(2z)\right)
+2k_0k_2\left(2k_0^2\partial_z^4+\partial_z^2\right),
\]
and $k_2$ is given by \eqref{e:3}. The solvability condition for this equation gives
\[
\nu_2 = \frac{c^2}{4X_2} >0
\]
which together with the fact that $\nu_{a,c}=O(a^2c^2)$ implies that $\nu_{a,c}>0$ and completes the proof of Theorem~\ref{t:1}.

\subsection{Proof of Theorem~\ref{t:3}}
\label{ss:pt3}

Consider $a$, $c$, and $\Lambda$ such that the results in Proposition~\ref{p:2} and Theorem~\ref{t:1} hold, and set $\lambda = k_{a,c}\Lambda$.

First, we claim that the operator $\lambda\partial_x-\mathcal C_{a,c}$ is a relatively compact perturbation of the asymptotic operator $\lambda\partial_x-\mathcal C_{a,c}^\infty$, when both operators act in $L^2(\R)$ with domains $H^6(\R)$. Indeed, the difference
\[
\mathcal G_{a,c} = \left(\lambda\partial_x-\mathcal C_{a,c}^\infty\right)
-\left(\lambda\partial_x-\mathcal C_{a,c}\right)
=
\left\{ \begin{array}{ll}
\partial_x^2(g_{a,c}^+\cdot), & \mbox{ for } x>0\\
\partial_x^2(g_{a,c}^-\cdot), & \mbox{ for } x<0
\end{array}\right. 
\] 
where
\[
g_{a,c}^\pm(x) =
u_{a,c}(x)-\varphi_{a,c}(x\pm \tau_{a,c})
\]
defines a closed operator in $L^2(\R)$ with domain $H^2(\R)$.
Since $g_{a,c}^\pm$ is a smooth function on $\R^\pm$ with 
$\lim_{x\to\pm \infty}\partial_x^j g_{a,c}^\pm(x) = 0$, $0\le j\le 4$, by Proposition~\ref{p:2}, using the compact embedding of $H^4(I)$ into $L^2(I)$ for any bounded interval $I$ and the continuity of $\mathcal G_{a,c}$ as an operator from $H^6(\R^\pm)$ to $H^4(\R^\pm)$, we conclude that for any bounded sequence $(f_n)_{n\geq1}\subset H^6(\R)$, the sequence 
$(\mathcal G_{a,c}f_n)_{n\geq1}\subset L^2(\R)$ contains a convergent subsequence. This implies that $\mathcal G_{a,c}$ is relatively compact with respect to $\lambda\partial_x-\mathcal C_{a,c}^\infty$ and proves the claim. As a consequence, the operators  $\lambda\partial_x-\mathcal C_{a,c}^\infty$ and  $\lambda\partial_x-\mathcal C_{a,c}$ have the same essential spectrum, so that it is enough to show that  $\lambda\partial_x-\mathcal C_{a,c}^\infty$ has negative essential spectrum.

According to Theorem~\ref{t:1}, there exists $\nu_*<0$ and a nontrivial $2\pi$-periodic smooth function $u_*$ such that 
\[
\left(\Lambda\partial_z-\mathcal B_{a,c}\right) u_* = \nu_*u_*.
\]
We set $v_*(x) = u_*(k_{a,c}x+\tau_{a,c})$, which solves the eigenvalue problem
\[
\left(\lambda\partial_x-\mathcal A_{a,c}^+\right) v_* = k_{a,c}^2\nu_*v_*,
\]
and consider  a cut-off function $\phi\in C_0^\infty(\R)$ such that
\[
\phi(x) = \left\{
\begin{array}{ll}
1, & \mbox{if } x\in[1,2]\\
0, & \mbox{if } x\in(-\infty,0]\cup [3,\infty)
\end{array}\right. .
\]
We define the sequence 
\[
v_n(x) = v_*(x)\phi_n(x),\quad n\geq1,
\]
where $\phi_n$ is the smooth function defined by
\[
\phi_n(x) = \left\{
\begin{array}{ll}
\phi(x), & \mbox{if } x\in[0,1]\\
1, & \mbox{if } x\in[1,n+1]\\
\phi(x-n+1), & \mbox{if } x\in[n+1,n+2]\\
0, & \mbox{if } x\in(-\infty,0]\cup [n+2,\infty)
\end{array}\right. .
\]
Since $v_*$ is a periodic function we have that  
$\|v_n\|\to\infty$ as $n\to\infty$, and for 
$\mathcal N_* = \lambda \partial_x -\mathcal C_{a,c}^\infty - k_{a,c}^2\nu_*$
we find
\[
\|\mathcal N_*v_n\|^2 = \int_0^1 |\mathcal N_*(v_*(x)\phi(x))|^2dx
+ \int_{n+1}^{n+2}  |\mathcal N_*(v_*(x)\phi(x-n+1))|^2dx \leq C_*,
\]
for any $n\geq1$, and some positive constant $C_*$ which does not depend on $n$. As a consequence, the operator $\mathcal N_*$ is not Fredholm, which implies that $k_{a,c}^2\nu_*<0$ belongs to the essential spectrum of $\lambda \partial_x -\mathcal C_{a,c}^\infty $. This completes the proof of Theorem~\ref{t:3}.

\appendix

\section{Existence of generalized solitary waves}
\label{a:exist_gw}

In this appendix we show how the results in Proposition~\ref{p:2} follow from the general result in \cite[Chapter 7, Theorem 7.1.18]{Lo}. In particular, this will also allow to recover the results in Proposition~\ref{p:1}.

We start by writing the equation \eqref{e:2} as a first order system
\begin{equation}\label{e:5}
\frac{dU}{dx} = \mathcal V(U,c),
\end{equation}
in which
\[
U=\begin{pmatrix}
u\\u_1\\u_2\\u_3
\end{pmatrix},\quad
\mathcal V(U,c)=\begin{pmatrix}
u_1\\u_2\\u_3\\-u_2+cu-\frac12 u^2
\end{pmatrix}.
\]
Notice that the system \eqref{e:5} is reversible, i.e., the vector field $\mathcal V$ anti-commutes with the reflection $S=\mathrm{diag}(1,-1,1,-1)$.

For any $c\in\R$, the system \eqref{e:5} possesses the equilibrium $U=0$. By linearizing at $U=0$ we find the Jacobian matrix
\[
J_c = \begin{pmatrix}
0&1&0&0\\
0&0&1&0\\
0&0&0&1\\
c&0&-1&0
\end{pmatrix},
\]
with eigenvalues $\nu$ satisfying
\[
\nu^4+\nu^2-c=0.
\]
At $c=0$, we find the double non-semi-simple eigenvalue $0$ and the simple eigenvalues $\pm\rmi$. We consider a basis $\{\varphi_0,\varphi_1,\varphi_+,\varphi_-\}$ consisting of eigenvectors and generalized eigenvectors of $J_0$, 
\[
\varphi_0 = \begin{pmatrix} 1\\0\\0\\0 \end{pmatrix},\quad
\varphi_1 = \begin{pmatrix} 0\\1\\0\\0 \end{pmatrix},\quad
\varphi_+ = \begin{pmatrix} -1\\-\rmi\\1\\ \rmi \end{pmatrix},\quad
\varphi_- = \begin{pmatrix} -1\\ \rmi\\1\\-\rmi \end{pmatrix},
\]
satisfying
\[
J_0\varphi_0=0,\quad J_0\varphi_1 = \varphi_0,\quad
J_0\varphi_+ = \rmi\varphi_+,\quad J_0\varphi_-=-\rmi\varphi_-,\quad
S\varphi_0=\varphi_0,
\]
together with the dual basis
\[
\varphi_0^* = \begin{pmatrix} 1\\0\\1\\0 \end{pmatrix},\quad
\varphi_1^* = \begin{pmatrix} 0\\1\\0\\1 \end{pmatrix},\quad
\varphi_+^* = \frac12\begin{pmatrix} 0\\0\\1\\ \rmi \end{pmatrix},\quad
\varphi_-^* = \frac12\begin{pmatrix} 0\\0\\1\\-\rmi \end{pmatrix}.
\]
Following \cite[Chapter 7]{Lo}, we compute the scalar products
\[
d_{10} = \langle D^2_{Uc}\mathcal V(0,0)\varphi_0,\varphi_1^* \rangle = 1,\quad
d_{20} = \langle D^2_{UU}\mathcal V(0,0)[\varphi_0,\varphi_0],\varphi_1^* \rangle =- 1,
\]
which are both non-zero. This shows that we are in the presence of a $0^2(\rmi)$ resonance and that the quadratic vector field is not degenerate.

The starting point in the construction of generalized solitary waves is a normal form transformation of \eqref{e:5}, followed by a scaling transformation \cite[Section 7.1.1]{Lo}. First, there exists a close to identity polynomial change of coordinates, analytically depending upon $c$, and preserving reversibility,
\[
U = \widetilde\alpha \varphi_0 +\widetilde\beta \varphi_1
+ \widetilde A (\Im\varphi_+) + \widetilde B (\Re\varphi_+) + \Phi(\widetilde Y,c), 
\]
in which $\widetilde Y = (\widetilde \alpha, \widetilde \beta, \widetilde A, \widetilde B)$ and $\Phi$ is a polynomial in $\widetilde Y$ with coefficients depending analytically upon $c$,
such that the system \eqref{e:5} is equivalent in a neighborhood of the origin to 
\[
\frac{d\widetilde Y}{dx} = \widetilde N(\widetilde Y,c) +
\widetilde R(\widetilde Y,c),
\]
where $\widetilde R(\widetilde Y,c) = O(|\widetilde Y|^3)$ and 
 $\widetilde N$ is the normal form of order $2$,
\[
\widetilde N(\widetilde Y,c) = 
\begin{pmatrix}
\widetilde \beta\\
d_1(c) c \widetilde \alpha + d_2(c) \widetilde \alpha^2 +
d_3(c) \left( \widetilde A^2 + \widetilde B^2\right)\\
-\widetilde B\left( 1 + c\omega_1(c) + m(c)\widetilde \alpha\right)\\
\widetilde A\left( 1+ c\omega_1(c) + m(c)\widetilde \alpha\right)
\end{pmatrix}.
\]
Here $d_1,d_2,d_3,\omega_1$ and $m$ are analytic functions of $c$ and
\[
d_1(0) = d_{10}=1,\quad d_2(0) = d_{20}=-1.
\]

Next, for $c>0$, we introduce the scaling
\[
x = c^{-1/2}y,\quad \widetilde \alpha = \frac32 c\alpha,\quad \widetilde \beta = \frac32 c^{3/2}\beta,\quad \widetilde A = cA,\quad \widetilde B = cB, 
\]
which leads to the system 
\begin{equation}\label{e:6}
\frac{dY}{dy} = N(Y,\sqrt c) + R(Y,\sqrt c),
\end{equation}
with
\[
N(Y,\sqrt c) = 
\begin{pmatrix}
\beta\\
\alpha -\frac32 \alpha^2 - \frac23
d_3(0) \left( A^2 +  B^2\right)\\
-B\left( \frac1{\sqrt c} + \omega_1(0)\sqrt c + m(0)\sqrt c \alpha\right)\\
A\left( \frac1{\sqrt c} + \omega_1(0)\sqrt c + m(0)\sqrt c \alpha\right)
\end{pmatrix},
\]
and $R(Y,\sqrt c)$ representing higher order terms.

The result in \cite[Theorem 7.1.4]{Lo} shows the existence of periodic orbits for \eqref{e:6}. Transforming back to the equation \eqref{e:2} this result is precisely the one stated in Proposition~\ref{p:1}, when restricting to positive values $c\in(0,c_0]$. Next, the result in \cite[Theorem 7.1.18]{Lo} shows the existence of reversible homoclinic orbits to small periodic orbits for the system \eqref{e:6}, provided the size of the periodic orbits is larger than an exponentially small critical size. For the equation \eqref{e:2} this leads to the result in Proposition~\ref{p:2}.

\noindent\\
{\bf Acknowledgements.} M. Haragus has been partially supported by the ANR project BoND (ANR-13-BS01-0009-01). E. Wahl\'{e}n was supported by the Swedish Research Council (grant no. 621-2012-3753).

\end{document}